\documentclass[12pt]{amsart}
\usepackage{fullpage,amsfonts,amsmath,amscd,amssymb}
\newtheorem*{thm}{Theorem}
\newtheorem*{prop}{Proposition}
\newtheorem*{lem}{Lemma}

\newtheorem*{rem}{Remark}

\newcommand{\hr}{\mathfrak{h}^{\text{reg}}}

\newcommand{\h}{{\mathfrak h}}

\newcommand{\ep}{\epsilon}
\newcommand{\irr}[1]{\textsf{Irrep}(#1)}

\newcommand{\Hom}{{\mbox{\text Hom}}}
\newcommand{\End}{{\mbox{\text End}}}

\newcommand{\hec}[1]{\mathcal{H}_W(#1)}

\DeclareMathOperator{\ed}{\End}

\title{On the quotient ring by diagonal invariants}
\author{Iain Gordon}
\address{Department of Mathematics, University of Glasgow, 
Glasgow, G12 8QW, U.K.} 
\email{ig@maths.gla.ac.uk}
\begin{document}
\begin{abstract} For a Weyl group, $W$, and its reflection representation $\h$, we find the character and Hilbert series for a quotient ring of $\mathbb{C}[\h\oplus \h^*]$ by an ideal containing the $W$--invariant polynomials without constant term. This confirms conjectures of Haiman.
\end{abstract}
\maketitle
\section{Introduction}
\subsection{} Let $W$ be a Weyl group, $\h$ its reflection representation, and $\mathbb{C}[\h]$ the ring of polynomial functions on $\h$. The action of $W$ on $\h$ extends to an action by algebra automorphisms on $\mathbb{C}[\h]$ by $w\cdot f(x) = f(w^{-1}\cdot x)$. It is a classical fact that the ring of invariants
$$\mathbb{C}[\h]^W = \{ f\in \mathbb{C}[\h] : w\cdot f = f \text{ for all $w\in W$}\}$$ 
is a polynomial ring in $\dim (\h)$ variables, \cite[Chapter 3]{hum}. Furthermore, the coinvariant ring
$$ \mathbb{C}[\h]^{co W} = \frac{\mathbb{C}[\h]}{\langle \mathbb{C}[\h]^W_+\rangle}$$ is a finite dimensional vector space, isomorphic as a $W$--module to the regular representation $\mathbb{C}W$ of $W$, \cite[Chapter 3]{hum}. (Here, $\langle \mathbb{C}[\h]^W_+\rangle$ denotes the ideal of $\mathbb{C}[\h]$ generated by invariant polynomials without constant term.) These important results are fundamental in algebraic geometry, the representation theory of finite simple groups of Lie type and the theory of Lie algebras.
\subsection{}
Recently, attention has focused on a ``double" analogue of the above results. The space $\h$ is replaced by $\h\oplus \h^*$, and its corresponding diagonal $W$--action. The orbit space $\h\oplus \h^*/W$ is a particularly interesting example of a symplectic singularity, of current interest in algebraic geometry, \cite{verb}. Moreover, the process of ``doubling up" is crucial to the study of Nakajima quiver varieties and preprojective algebras, the Hilbert scheme of points in the plane, \cite{nak}, and in Lie theoretic work on principal nilpotent pairs, \cite{ginz}.
\subsection{}
The ring of invariants $\mathbb{C}[\h\oplus \h^*]^W$ is never smooth and for general $W$, its generators and relations are poorly understood, \cite{wallach}. It is, however, expected that the ring of coinvariants $$\mathbb{C}[\h\oplus \h^*]^{co W}=\frac{\mathbb{C}[\h\oplus \h^*]}{\langle \mathbb{C}[\h\oplus \h^*]_+^W\rangle }$$ should display interesting combinatorial properties, relating to its $W$--action and its natural grading arising from setting $\deg (x)=1$ and $\deg (y) = -1$ for $x\in \h^*$ and $y\in \h$, \cite{haiman}. 

In the case $W=\mathfrak{S}_n$, the work of Haiman on the $n!$--conjecture shows that the Hilbert series of $\mathbb{C}[\h\oplus \h^*]^{co W}$ is $t^{-n(n-1)/2}(1+t +\cdots + t^n)^{n-1}$, and determines its decomposition as an $\mathfrak{S}_n$--module, \cite{haiman2}. This confirms several conjectures in \cite{haiman}.
\subsection{}
In this paper we extend the above results to all Weyl groups. Our main theorem is the following, confirming a conjecture of Haiman \cite[Section 7]{haiman}.
\begin{thm}
\label{maintheorem}
Let $W$ be a Weyl group. Let $n$ be the rank of $W$ and $h$ the Coxeter number. Let $D_W= \mathbb{C}[\h\oplus \h^*]^{co W}$. Then there exists a $W$--stable quotient ring $R_W$ of $D_W$ satisfying the following properties:
\begin{enumerate} 
\item $\dim R_W = (h+1)^n$;
\item $R_W$ is $\mathbb{Z}$--graded with Hilbert series $t^{-N}(1+t+\cdots +t^h)^n$;

\item The image of $\mathbb{C}[\h]$ in $R_W$ is the classical coinvariant algebra, $\mathbb{C}[\h]^{co W}$;
\item As a $W$--module $R_W\otimes \ep$ is isomorphic to the permutation representation of $W$ on the reduction of the root lattice modulo $h+1$, written $Q/(h+1)Q$.
\end{enumerate}
\end{thm} 
It is known that for $W$ not of type $A$, $R_W$ is usually a \textit{proper} quotient of $D_W$. For dihedral groups the theorem is known by calculations of Alfano and Reiner. 

\subsection{}
Let us make some comments on the method of proof, which is rather indirect. Our strategy is to find a non--commutative deformation of $R_W$ which has all the properties discussed in the theorem. This non--commutative deformation is a simple module of a rational Cherednik algebra associated to $W$. The rational Cherednik algebras, introduced by Etingof and Ginzburg \cite{EG}, are members of a certain flat family of deformations of the skew group algebra $\mathbb{C}[\h\oplus \h^*]\ast W$. The parameter set for this family is the $W$--invariant functions $c: R\longrightarrow\mathbb{C}$, where $R$ denotes the set of roots attached to $W$. We denote the  corresponding algebra by $H_c$. It is known that $\mathbb{C}W \subset H_c$ for all $c$.

Of crucial importance is the existence of a faithful ``Dunkl representation" for these algebras, allowing us  to relate the representation theory of a rational Cherednik algebra $H_c$ through monodromy to the representation theory of the Hecke algebra associated to $W$ at parameter $q= exp(2\pi ic)$. 

Using known results on the representation theory of Hecke algebra at roots of unity, we can describe quite well a simple module, $L$, for $H_c$ at a particular value of $c$.  In particular, it is possible to understand the decomposition of $L$ as a $W$--module, and describe a certain $\mathbb{Z}$--grading on it. 

Let $e_{\ep}\in \mathbb{C}W$ be the sign idempotent. Define the \textit{spherical algebra} to be $e_{\ep}H_ce_{\ep}$ for any $c$ and note that $H_ce_{\ep}$ is a $(H_c, e_{\ep}H_ce_{\ep})$--bimodule. The spherical algebra has a one--dimensional module $\mathbb{C}$. It turns out that $L$ has an alternative description as the induced module $H_ce_{\ep} \otimes_{e_{\ep}H_ce_{\ep}} \mathbb{C}$. Passing to the associated graded module we find a surjective homomorphism from $\mathbb{C}[\h\oplus \h^*]e_{\ep}\otimes_{\mathbb{C}[\h\oplus \h^*]^W e_{\ep}} \mathbb{C}$ to $gr L$. For Theorem \ref{maintheorem} we set $R_W = gr L$.

\subsection{} It seems likely that the proof can be adapted to deal with a larger class of complex reflection groups. In particular for the wreath products $\mathbb{Z}_m\wr \mathfrak{S}_n$, that is complex reflection groups of type $G(m,1,n)$, results of Weyl on the generation of $\mathbb{C}[\h\oplus \h^*]^W$ by generalised polarisations should allow one to prove an analogue of the important Proposition \ref{shift}. Using the work of Dunkl and Opdam, \cite[Section 3.5]{do} on shifts of Dunkl operators, and the work of Ariki, \cite{ariki}, and Brou\'{e}, Malle and Rouquier, \cite{bmr}, on cyclotomic Hecke algebras and their representations, it should be possible to extend Theorem \ref{maintheorem} to this larger class of groups. 

\subsection{} The paper is organised as follows. In Section 2 we recall basic definitions from rational Cherednik algebras. Section 3 discusses Hecke algebras and their relation to rational Cherednik algebras. In Section 4 we find the finite dimensional simple $H_c$--module described above. Finally, in Section 5, we relate the earlier sections to the coinvariants of $W$ on $\mathbb{C}[\h\oplus \h^*]$, proving Theorem \ref{maintheorem}. 

\subsection{}  In recent work, \cite{BEG2}, Berest, Etingof and Ginzburg have found a character formula for \textit{all} finite dimensional simple $H_c$--modules in type $A$, and for many in other types. Their techniques, which were discovered independently, are similar to those used here and point towards links with Hilbert schemes on the plane and other ``doubled'' geometric objects. We thank Ginzburg for bringing this work to our attention. Moreover, we are indebted to Etingof for showing how to extend our earlier version of Proposition \ref{shift} and of Theorem \ref{maintheorem}(4) from classical Weyl groups to the exceptional types. Without this, our results would be less complete.

\section{Rational Cherednik algebras}
\subsection{} The following is due to Etingof and Ginzburg, \cite{EG}. Let $W$ be a Weyl group, and $\mathfrak{h}$ its reflection representation 
over $\mathbb{C}$. Let $R\subset \h^*$ be the root system corresponding to $W$. To each $W$--invariant function $c:R\longrightarrow \mathbb{C}$ one can attach an associative algebra $H_c$, called the \textit{rational Cherednik algebra}. Given $\alpha\in R$, write $\alpha^{\vee}\in \h$ for the corresponding coroot and $s_{\alpha}\in GL(\h)$ for the reflection corresponding to $\alpha$. The rational Cherednik algebra $H_c$ is generated by the vector spaces $\h, \h^*$, and the group $W$, with defining relations given by
\begin{eqnarray*} w\cdot x\cdot w^{-1} = w(x), \quad w\cdot y\cdot w^{-1} = w(y), \qquad &&\text{for all }y\in \h, x\in \h^*, w\in W \\ x_1\cdot x_2 = x_2\cdot x_1, \quad y_1\cdot y_2 = y_2\cdot y_1, \qquad && \text{for all }y_1,y_2\in \h, x_1,x_2\in \h^* \\ y\cdot x - x\cdot y = \langle y,x \rangle - \frac{1}{2}\sum_{\alpha\in R} c_{\alpha} \cdot \langle y,\alpha \rangle\langle \alpha^{\vee}, x\rangle \cdot s_{\alpha}, \qquad && \text{for all }y\in \h, x\in \h^*.\end{eqnarray*}
\subsection{} The elements $x\in\h^*$ generate a subalgebra $\mathbb{C}[\h]\subset H_c$ of polynomial functions on $\h$, the elements $y\in \h$ generate a subalgebra $\mathbb{C}[\h^*]\subset H_c$, and the elements $w\in W$ span a copy of the group algebra $\mathbb{C}W$ sitting naturally inside $H_c$. By \cite[Theorem 1.3]{EG} there is a Poincar\'{e}--Birkhoff--Witt isomorphism \begin{equation} \label{PBW} \begin{CD} \mathbb{C}[\h]\otimes \mathbb{C}W \otimes \mathbb{C}[\h^*] @>\sim >> H_c.\end{CD}\end{equation}

\subsection{} Let $\{ x_i\}$ and $\{ y_i \}$ be a pair of dual bases of $\h^*$ and $\h$ respectively. Viewing $\h$ and $\h^*$ as subspaces of $H_c$, we let $\mathbf{h} = \frac{1}{2}\sum_i (x_iy_i + y_ix_i)\in H_c$, which is independent of the choice of bases. Let $( \, , \, )$ be the complex bilinear form on $\h$ extending the Euclidean inner product. Define $\mathbf{x}^2\in \mathbb{C}[\h]$ to be the squared norm function $x\mapsto (x,x)$ and let $\mathbf{y}^2\in \mathbb{C}[\h^*]$ be defined similarly. By \cite[(2.6)]{BEG}, $\{\mathbf{x}^2, \mathbf{h}, \mathbf{y}^2\}$ forms an $\mathfrak{sl}_2$--triple in the algebra $H_c$.

\subsection{A filtration}\label{filt} There exists a filtration on $H_c$ with $\deg(x) = \deg(y) = 1$ and $\deg(w)=0$ for $x\in\h^*$, $y\in \h$ and $w\in W$. By the Poincar\'{e}-Birkhoff-Witt isomorphism, the associated graded ring of $H_c$ with respect to this filtration, $gr H_c$, equals $\mathbb{C}[\h\oplus \h^*]\ast W$. Using the non--commutativity of $H_c$, a standard procedure, \cite[Section 2]{EG}, produces a Poisson bracket on $\mathbb{C}[\h\oplus \h^*]^W$. This agrees with the Poisson bracket induced by the canonical $W$--invariant symplectic form on $\h\oplus \h^*$, \cite[Lemma 2.23]{EG}.

\subsection{Category $\mathcal{O}_c$} 
We recall the definition of a special subcategory of $H_c$--mod from \cite[Definition 2.4]{BEG}. Let $\mathcal{O}_c$ be the abelian category of finitely--generated $H_c$--modules $M$, such that the action on $M$ of the subalgebra $\mathbb{C}[\h^*]\subset H_c$ is locally finite, that is $\dim \mathbb{C}[\h^*]\cdot m< \infty$ for any $m\in M$, and such that for any $P\in \mathbb{C}[\h^*]^W$ the action of $P- P(0)$ is locally nilpotent. Let $\irr{W}$ denote the set of simple $W$--modules, up to isomorphism. Given $\tau \in \irr{W}$, we define an object of $\mathcal{O}_c$, called the \textit{standard} module $M_c(\tau)$, to be the induced module $$M_c(\tau) = H_c\otimes_{\mathbb{C}[\h^*]\ast W} \tau,$$ where $\mathbb{C}[\h^*]\ast W$ acts on $\tau$ by sending $pw\cdot m = p(0) (w\cdot m)$, for $p\in \mathbb{C}[\h^*]$, $w\in W$ and $m\in \tau$. It is shown in \cite[Section 2]{BEG} that each $M_c(\tau)$ has a unique simple quotient $L_c(\tau)$, that the set $\{ L_c(\tau) : \tau \in \irr{W}\}$ provides a complete list of non--isomorphic simple objects in $\mathcal{O}_c$, and that every object in $\mathcal{O}_c$ has finite length. Note that it follows from the Poincar\'{e}--Birkhoff--Witt theorem that the standard module $M_c(\tau)$ is a free $\mathbb{C}[\h]$-module of rank $\dim(\tau)$.

\section{Hecke Algebras}
\subsection{}
Let $B_W$ be the braid group of $W$, in other words the fundamental group of the variety $\hr/W$. Fix a point $\ast \in \hr$ and for each simple reflection $s_{\alpha}\in W$, let $T_{\alpha}$ be the class in $B_W$ corresponding to a straight path from the point $\ast$ to the point $s_{\alpha}(\ast)$ with an inserted anti--clockwise semi--circle around the hyperplane $\alpha = 0$. Given $\mathbf{q} : R\longrightarrow \mathbb{C}^{\ast}$ ($\alpha\mapsto q_{\alpha}$), a $W$--invariant function, we define $\hec{\mathbf{q}}$ to be the quotient of the group algebra $\mathbb{C}[B_W]$ by the ideal generated by $(T_{\alpha}-1)(T_{\alpha} - q_{\alpha})$. The algebra $\hec{\mathbf{q}}$ is the \textit{Hecke algebra} of $W$. 
\subsection{}
It is well--known that $\dim \hec{\mathbf{q}} = |W|$. When $\mathbf{q}=\mathbf{1}$, $\hec{\mathbf{q}} = \mathbb{C}W$. For a generic choice of $\mathbf{q}$ the algebra $\hec{\mathbf{q}}$ is semisimple and isomorphic to $\mathbb{C}W$. 
For any choice of $\mathbf{q}$, there exist a set of $\hec{\mathbf{q}}$-modules labelled by the simple $W$--modules, \cite{mathas}. These are called \textit{Specht modules} and written $S_{\mathbf{q}}(\lambda)$ for $\lambda \in \irr{W}$. These are simple--headed and give a complete set of non--isomorphic irreducible $\hec{\mathbf{q}}$--modules in the semisimple case. If $\hec{\mathbf{q}}$ is not semisimple the composition factors of $S_{\mathbf{q}}(\lambda)$ are the entries of the decomposition matrix for $\hec{\mathbf{q}}$.  
\subsection{}
We need to recall a few results about Dunkl operators, which link rational Cherednik algebras to Hecke algebras, see \cite[Section 2]{BEG}. According to \cite[Proposition 4.5]{EG}, the algebra $H_c$ has a faithful ``Dunkl representation", an injective
algebra homomorphism $H_c \longrightarrow D(\hr )\ast W$. This allows us to consider  $M|_{\hr}$, the localisation of an $H_c$--module $M$, as a $W$--equivariant $D$--module on $\hr$. 

In particular, the standard module $M_c(\tau)|_{\hr}$, viewed as a $D$-module on $\hr$, is the trivial vector bundle $\mathbb{C}[\hr]\otimes \tau$ equipped with a flat connection. It is well--known that the connection is the Knizhnik-Zamolodchikov connection with values in $\tau$. Moreover, the monodromy representation of the fundamental group $B_W=\pi_1(\hr/W)$ corresponding to this connection factors through the Hecke algebra $\hec{exp(2\pi ic)}$.

\subsection{}
\label{orlem}
In general, given a $W$--equivariant vector bundle $M$ on $\hr$ with a flat connection, the germs of the horizontal holomorphic sections of $M$ form a locally constant sheaf on $\hr/W$. Let $Mon(M)$ be the corresponding monodromy representation of the fundamental group $\pi_1(\hr/W)$ in the fibre over $\ast$, where $\ast$ is some fixed point in $\hr/W$. The following Lemma was proved by Opdam--Rouquier and presented in \cite[Lemma 2.10]{BEG}.
\begin{lem}
Let $N$ be an object of $\mathcal{O}_c$ which is torsion--free over the subalgebra $\mathbb{C}[\h]\subset H_c$. Then for any $M$ belonging to $\mathcal{O}_c$, the canonical map $$\Hom_{H_c}(M,N) \longrightarrow \Hom_{B_W}(Mon(M|_{\hr}),Mon(N|_{\hr}))$$ is injective.
\end{lem}
In particular, this lemma applies to the case $N= M_c(\tau)$ for any $\tau\in\irr{W}$.

\section{A simple module for rational Cherednik algebras}

\subsection{}
Let $\ep$ be the sign representation of $W$ and $e_{\ep} = |W|^{-1}\sum_w \ep(w)w$ the sign idempotent. Let $e= |W|^{-1}\sum_w w$, the trivial idempotent. Let $h$ be the Coxeter number of $W$. Given $0\leq k\leq n$, set $\h_k=\wedge^k\h$, the $k$--th exterior power of $\h$. In particular $\h_1 = \h$ and $\h_0 = \textsf{triv}$. Thanks to \cite[Theorem 9.13]{curkil} $\h_k\in \irr{W}$. Let $n=\dim \h$. Recall $R$ is the set of roots associated with $W$, and let $R_+$ be a set of positive roots. Set $N=|R_+|$.

\subsection{}
\label{centchar}
Each module $M$ in category $\mathcal{O}_c$ splits into a direct sum of generalised $\mathbf{h}$--eigenspaces, \cite[Section 2]{guay}.  Moreover, the action of $\mathbf{h}$ on the standard and simple modules in $\mathcal{O}_c$ is diagonalisable. To describe the $\mathbf{h}$--eigenspaces of $M_c(\tau)$ we introduce some notation. Let $$\kappa_c = \sum_{\alpha\in R_+} c_{\alpha}(1-s_{\alpha})$$ be the central element of $\mathbb{C}W$ canonically attached to $c\in\mathbb{C}[R]^W$. This element acts by a scalar, say $\kappa_c(\tau)$, on each irreducible representation $\tau\in \irr{W}$. Thanks to \cite[Section 2]{guay}, for $p\in\mathbb{C}[\h]_m$ and $v\in \tau$ we have in $M_c(\tau)$ \begin{equation} \label{weechar} \mathbf{h}(p\otimes v) = (m+\frac{n}{2} +\kappa_c(\tau) - \sum_{\alpha\in R_+} c_{\alpha})p\otimes v.\end{equation}
\begin{lem} 
Keep the above notation. If $c$ is constant, then $\kappa_c(\h_k) = hck$.
\end{lem}
\begin{proof}
Recall $N=|R_+|$ and $\h_k = \wedge^k \h$. Since $\kappa_c$ is central, for all $\tau\in\irr{W}$ we have \begin{equation}\label{char} \kappa_c(\tau) = \frac{\chi_{\tau}(\kappa_c)}{\dim \tau} = Nc - \frac{c\sum \chi_{\tau}(s_{\alpha})}{\dim \tau},\end{equation} where $\chi_{\tau}$ is the character of $\tau$. 

For any $\alpha\in R$, the reflection $s_{\alpha}$ fixes a hyperplane pointwise and acts as $-1$ on its orthogonal complement, so that $\h = F\oplus S$, where $\dim F=n-1$ and $\dim S=1$. Thus, the action of $s_{\alpha}$ on $\h_k$ is described by $\h_k = \wedge^k F \oplus \wedge^{k-1}F\otimes S$. We deduce that $$\chi_{\h_k}(s_{\alpha}) = \binom{n-1}{k} - \binom{n-1}{k-1}.$$ Combining this with (\ref{char}) yields $$\kappa_c(\h_k) = Nc - \frac{Nc\left(\binom{n-1}{k} - \binom{n-1}{k-1}\right)}{\binom{n}{k}} = Nc - \frac{Nc(n-2k)}{n} = \frac{2Nck}{n}.$$ Since the Coxeter number of $W$ equals $2N/n$, \cite[Chapter 3]{hum}, the lemma follows.
\end{proof}

\subsection{}
In the following lemma, we make use of the relationship between rational Cherednik algebras and Hecke algebras. 
\begin{lem}
\label{maps}
Suppose $c=(1+mh)/h$ for $m$ a positive integer. Let $\tau, \mu\in \irr{W}$ be non--isomorphic and suppose that there is a non--zero $H_c$--homomorphism $M_c(\tau)\longrightarrow M_c(\mu)$. Then $\tau \cong \h_i$ and $\mu\cong \h_j$ for some non--negative integers $i$ and $j$.
\end{lem}
\begin{proof}
By \cite[Section 6]{DJO} the composition multiplicities of the $\hec{q}$--modules $Mon(M_c(\mu)|_{\hr})$ are described by the decomposition matrix for $\hec{q}$. The decomposition matrix for $\mathbf{q}= \exp(2\pi i c)$ is described in \cite[Section 6.4]{gecketal}. In particular, \cite[Lemma 6.5]{gecketal} shows that the Specht modules associated to $\lambda$, for $\lambda \neq \h_i$ ($0\leq i\leq n$), are irreducible and projective. Hence, to each such $\lambda$, there corresponds $\mu$ such that $Mon(M_c(\mu)|_{\hr})$  is irreducible and projective. By \cite[Corollary 3.5]{dunkl} we must have that $\mu = \lambda$. Thus we have shown that $Mon(M_c(\lambda)|_{\hr})$ is an irreducible and projective $\hec{q}$--module for $\lambda \neq \h_i$, $0\leq i\leq n$. Furthermore, by \cite[Theorem 6.6]{gecketal} there is a unique non--trivial block for $\hec{q}$ containing the Specht modules associated to $\h_i$, so we deduce that $Mon(M_c(\h_i)|_{\hr})$ is not isomorphic to $Mon(M_c(\lambda)|_{\hr})$.

We have $\dim \Hom_{H_c}(M_c(\tau),M_c(\mu)) \leq \Hom_{\hec{q}}(Mon(M_c(\tau)|_{\hr}), Mon(M_c(\mu)|_{\hr}))$ by Lemma \ref{orlem}. The result follows from the above paragraph.
\end{proof}

\subsection{}
\label{ppsign}
Let $e_1,\ldots ,e_n$ be the exponents of $W$ and $d_i=e_i+1$, the degrees of the fundamental invariants for $W$, \cite[Chapter 3]{hum}. Let $$p= \prod_{i=1}^n (1-t^{d_i})^{-1},$$ the Hilbert series of the invariant ring $\mathbb{C}[\h]^W$.
\begin{lem}
Suppose $c=(1+mh)/h$ for $m$ a positive integer. The Hilbert series of $e_{\ep}M_c(\h_{n-i})$ with respect to the $\mathbf{h}$--eigenspaces is 
$$ p(e_{\ep}M_c(\h_{n-i}), t) = t^{-mN +(n-i)(mh+1)}p \sum_{1\leq j_1<\cdots <j_i \leq n} t^{e_{j_1}+\cdots +e_{j_i}}.$$
\end{lem}
\begin{proof}
As a $\mathbb{C}[\h]\ast W$--module we have $M_c(\h_{n-i}) \cong \mathbb{C}[\h]\otimes \h_{n-i}$. Hence, the sign module appears in $M_c(\h_{n-i})$ with multiplicity one for each appearance of $(\h_{n-i})^* \otimes \ep$ in $\mathbb{C}[\h]$. Since $\h_n = \ep$, there is a homomorphism $\h_{n-i} \wedge \h_i \longrightarrow \ep$, and hence an isomorphism $\h_i \cong (\h_{n-i})^* \otimes \ep$. 

To determine the appearances of $\h_i$ in $\mathbb{C}[\h]$, we study $$\mathbb{C}[\h]^{co W} = \frac{\mathbb{C}[\h]}{\langle \mathbb{C}[\h]_+^W\rangle}.$$ By \cite{sol} there is a copy of $\h$ in degree $e_i$ for $1\leq i\leq n$. Since $\mathbb{C}[\h]^{co W}$ is isomorphic to the regular representation of $W$, we know that $\h_i$ appears $\binom{n}{i}$ times. Taking $i$ distinct copies of $\h$ in degrees $e_{j_1},\ldots ,e_{j_i}$ we obtain a copy of $\h_i$ in degree $e_{j_1}+\cdots +e_{j_i}$ by taking their wedge product. Running through all such products yields $\binom{n}{i}$ copies of $\h_i$. Since the wedge product of the $n$ distinct copies of $\h$ above yields the unique copy of $\ep$ in $\mathbb{C}[\h]^{co W}$, all the copies of $\h_i$ described above are distinct. Hence the Hilbert series of $\h_i$ in $\mathbb{C}[\h]^{co W}$ is given by $ \sum_{1\leq j_1<\cdots <j_i \leq n} t^{e_{j_1}+\cdots +e_{j_i}}$. 

It follows from the $W$--invariant graded decomposition $$\mathbb{C}[\h] \cong \mathbb{C}[\h]^{co W}\otimes \mathbb{C}[\h]^W$$ that the Hilbert series for $\h_i$ in $\mathbb{C}[\h]$ is given by $p\sum_{1\leq j_1<\cdots <j_i \leq n} t^{e_{j_1}+\cdots +e_{j_i}}.$ 

The lowest $\mathbf{h}$--eigenspace in $M_c(\h_{n-i})$ can be calculated from formula (\ref{weechar}) on $1\otimes v$. Since $c=(1+mh)/h$ and $nh/2 = |R_+|$, Lemma \ref{centchar} shows that this equals $-mN + (n-i)(1+mh).$ The lemma follows.
\end{proof}

\subsection{}
\label{faithful}
The following lemma is crucial.
\begin{lem}
Suppose $c=(1+mh)/h$ for $m$ a positive integer. Then $e_{\ep}L_c(\h_i) \neq 0$.
\end{lem}
\begin{proof}
If $L_c(\h_i)=M_c(\h_i)$ the claim is clear, so suppose $M_c(\h_i)$ is not simple. Let $R_c(\h_i)$ be its radical (unique maximal submodule). Then $M_c(\h_i)/R_c(\h_i) = L_c(\h_i)$. Moreover, the lowest $\mathbf{h}$--eigenspace of $R_c(\h_i)$ yields a homomorphism $M_c(\tau)\longrightarrow R_c(\h_i)$ for some $\tau\in\irr{W}$. By Lemma \ref{maps} we know that $\tau$ must have the form $\h_j$ for some $j$. By \cite[Section 3]{guay} we have $\kappa_c(\h_i) < \kappa_c(\h_j)$ which by Lemma \ref{centchar} shows that $j> i$. By Lemma \ref{ppsign}, the difference in the lowest $\mathbf{h}$--eigenspaces of $M_c(\h_i)$ and $M_c(\h_j)$ which contain a copy of $\ep$ is at least $mh+1 - e_i$ for some $i$. As all exponents $e_i$ are less than the Coxeter number, $M_c(\h_j)$ cannot map onto the lowest $\mathbf{h}$--eigenspace of $M_c(\h_i)$ containing a copy of $\ep$. Thus means that $R_c(\h_i)$ cannot contain all copies of $\ep$.  
\end{proof}

\subsection{}
\label{oned}
In this section we consider only the value $c=1/h$.
\begin{lem} Let $c=1/h$. Then $H_c$ has a one dimensional module, whose restriction to $W$ is the trivial module. 
\end{lem}
\begin{proof}
For type $A$ this follows from \cite[Remark following Proposition 5.7]{BEG}. We will prove this for type $B$, the other cases being similar.

Let $\{x_1, \ldots ,x_n\}$ be the standard basis for $\h^*$ and $\{y_1,\ldots ,y_n\}$ for $\h$. Recall, \cite[Chapter 2]{hum}, that $W=\mathbb{Z}_2\wr \mathfrak{S}_n$ and that the positive roots of $W$ can be chosen to be the elements of $\h$ of the form $x_i\pm x_j$ for $1\leq i < j\leq n$ (short roots) and $2x_i$ for $1\leq i\leq n$ (long roots). 

Let $c_s$ and $c_l$ be the values of $c:R \longrightarrow \mathbb{C}$ on the short and long roots respectively. Calculation gives the following formula:
\begin{equation*}
[y_i, x_j] = \delta_{ij} + \begin{cases} -c_s\sum_{t\neq i}(s_{x_{i+t}} + s_{x_{i-t}}) - 2c_ls_{x_i} \quad &\text{if $i=j$}, \\ -c_s(s_{x_{i+j}}-s_{x_{i-j}}) & \text{if $i\neq j$.} \end{cases}  
\end{equation*}
The trivial representation sends the group elements to $1$, so setting $c_s = c_l =\frac{1}{2n}=1/h$ proves that $[y_i,x_j]\mapsto 0$. The lemma follows.
\end{proof}
\begin{rem} For $c=1/h$ it can be shown using a Koszul resolution that we have $$[M_c(\h_i): L_c(\h_j)] = \begin{cases} 1 &\quad \text{if $j=i,i+1$}\\ 0&\quad \text{otherwise,}\end{cases} $$ and $M_c(\tau) = L_c(\tau)$ if $\tau\neq \h_i$ for all $i$. This is greatly generalised in \cite{BEG2}.
\end{rem}

\subsection{} 
\label{altsum}
Consider the functor $F: H_c-mod \longrightarrow e_{\ep}H_ce_{\ep}-mod$ which sends $M$ to $e_{\ep}M$. Since $e_{\ep}$ is an idempotent, this functor is exact. All objects in category $\mathcal{O}_c$ have finite length, hence Lemma \ref{faithful} shows that if $c=(1+mh)/h$ then $M_c(\tau)$ is simple for $\tau\neq \h_i$ and $[M_c(\h_i) : L_c(\h_j)] = [e_{\ep}M_c(\h_i): e_{\ep}L_c(\h_j)]$. Let $D$ be the $n+1\times n+1$ matrix whose $i,j$--th entry is given by $[M_c(\h_i): L_c(\h_j)]$. By Lemma \ref{centchar} and \cite[Section 3, Proof of Theorem 8]{guay} this is an upper triangular matrix. This means that there is a unique way to write $e_{\ep}L_c(\h_j)$ in terms of $e_{\ep}M_c(\h_i)$.
\begin{lem}
Suppose that $c= (1+h)/h$. We have $\sum_{i=0}^n (-1)^i p(e_{\ep}(M_c(\h_i),t) = 1$.
\end{lem}
\begin{proof}
By Lemma \ref{ppsign} we have \begin{eqnarray*} \sum_{i=0}^n (-1)^i p(e_{\ep}(M_c(\h_i)) &=& p\sum_{i=0}^n (-1)^i t^{-N +i(h+1)}\sum_{1\leq j_1< \cdots < j_{n-i} \leq n} t^{e_{j_1}+\cdots +e_{j_{n-i}}} \\ &= & p t^{-N}\prod_{k=1}^n(t^{e_{k}} - t^{h+1}) \\ &=& pt^{-N}\prod_{k=1}^nt^{e_{k}}(1-t^{h+1-e_{k}}).\end{eqnarray*}
By \cite[Chapter 3]{hum} the sum $\sum_{k=1}^n e_k = N$, whilst the multiset $\{ h+1-e_k : 1\leq k\leq n\}$ equals $\{ d_k: 1\leq k\leq n\}$. We deduce that \begin{eqnarray*} \sum_{i=0}^n (-1)^n p(e_{\ep}(M_c(\h_i)) = p\prod_{k=1}^n (1-t^{d_k}) = 1,\end{eqnarray*} as required.
\end{proof}

\subsection{Shifting}
\label{shift}
We now need an analogue of a theorem of Berest, Etingof and Ginzburg for non--regular values of $c$. Let $\mathbf{1}:R \longrightarrow \mathbb{C}$ be the map taking the value $1$ everywhere.  To indicate provenance, we give two versions of the following Proposition. The second is due to Berest, Etingof and Ginzburg and will appear in \cite{BEG2}. We thank the authors for allowing us to reproduce the result here.

\begin{prop}[Version 1] Let $W$ be a Weyl group of type $A,B,D$ or $G_2$. Then for any function $c$ there is an algebra isomorphism $e_{\ep}H_ce_{\ep} \cong eH_{c-\mathbf{1}}e$.
\end{prop}
\begin{proof} 
The proof follows the same lines as \cite[Proposition 4.11]{BEG} with one exception. We need to know for \textit{all} values of $c$ that $eH_ce$ is generated by its positive part $\mathbb{C}[\h]^We$ and the element $\mathbf{y}^2$. This follows from \cite[Appendix 2]{wallach}  which shows that diagonal invariants of $W$ on $\mathbb{C}[\h\oplus \h^*]$ are generated by $\mathbb{C}[\h]^W$ and $\mathbb{C}[\h^*]^W$ and their Poisson brackets. Since $gr (eH_ce) = eS(\h \oplus \h^*)e \cong S(\h \oplus \h^*)^W$ and the Poisson bracket is induced by the commutator in $eH_ce$, Section \ref{filt}, it follows that $eH_ce$ is generated by $\mathbb{C}[\h]^We$ and $\mathbb{C}[\h^*]^We$. Finally, \cite[Corollary 4.9]{BEG} shows that $eH_ce$ is generated by $\mathbb{C}[\h^*]^We$ and $\mathbf{y}^2e$.
\end{proof}  
\begin{prop}[Version 2, \cite{BEG2}] For any function $c$ there is an algebra isomorphism $e_{\ep}H_ce_{\ep} \cong eH_{c-\mathbf{1}}e$.
\end{prop}
\begin{proof}
Let $\theta_c : eH_{c-\mathbf{1}}e \longrightarrow \ed(\mathbb{C}[\h])$ be the spherical Harish--Chandra homomorphism, and $\theta_c^-: e_{\ep}H_ce_{\ep} \longrightarrow \ed(\mathbb{C}[\h])$ the antispherical Harish--Chandra homomorphism, \cite[Section 4]{EG}. The image of these maps lies in differential operators on $\hr$, preserving $\mathbb{C}[\h]$. Give $H_c$ the filtration induced by setting $\deg (x) = \deg (w) = 0$ and $\deg (y) =1$ for $x\in \h^*, y\in \h$ and $w\in W$, and put the filtration on differential operators induced by order. By \cite[Proposition 4.5]{EG} both $\theta_c$ and $\theta_c^-$ are flat families of injective filtration preserving homomorphisms, such that the associated graded maps are injective. Since, by \cite[Propositions 4.10 and 4.11]{BEG}, the images of $\theta_c$ and $\theta_c^-$ are equal for generic values of $c$, the images are equal for all values of $c$. The isomorphism from $e_{\ep}H_ce_{\ep}$ to $eH_{c-\mathbf{1}}e$ is thus given by $\theta_c^{-1}\theta_c^-$.
\end{proof}
\subsection{}
\label{hs}
For $c=(1+h)/h$ we can now describe the Hilbert series of $L_c(\textsf{triv})$.
\begin{thm}
Set $c=(1+h)/h$. The simple $H_c$--module $L_c(\textsf{triv})$ is finite dimensional. Its Hilbert series is given by $$p(L_c(\textsf{triv}),c) = t^{-N}(1+t+t^2+\cdots + t^{h})^n.$$
\end{thm}  
\begin{proof}
Under the shift isomorphism of Proposition \ref{shift}(Version 2) $\mathbf{x}^2e$ is sent to $\mathbf{x}^2e_{\ep}$ and $\mathbf{y}^2e$ to $\mathbf{y}^2e_{\ep}$, \cite[Proof of Proposition 4.11]{BEG}. Since $\{ \mathbf{x}^2e, \mathbf{h}e, \mathbf{y}^2e \}$ and $\{ \mathbf{x}^2e_{\ep} , \mathbf{h}e_{\ep}, \mathbf{y}^2e_{\ep}\}$ form $\mathfrak{sl}_2$--triples in $eH_{c-\mathbf{1}}e$ and $e_{\ep}H_ce_{\ep}$ respectively, it follows that the shift isomorphism sends $\mathbf{h}e= [ \mathbf{x}^2e, \mathbf{y}^2e]$ to $\mathbf{h}e_{\ep}=[ \mathbf{x}^2e_{\ep}, \mathbf{y}^2e_{\ep}]$. 

Thanks to Lemma \ref{oned} and Proposition \ref{shift}(Version 2) there is a one--dimensional $e_{\ep}H_ce_{\ep}$--module. Call this one--dimensional module $\mathbb{C}$. By the above paragraph $\mathbf{h}e_{\ep}$ acts on $\mathbb{C}$ with weight zero since $\mathbf{h}e$ does so on the trivial module. Thus Lemma \ref{altsum} shows how to write $\mathbb{C}$ in terms of the standard modules.  

Define a functor $G: e_{\ep}H_ce_{\ep}-mod \longrightarrow H_c-mod$ sending $M$ to $G(M) =  H_ce_{\ep} \otimes_{e_{\ep}H_ce_{\ep}} M$. Since $H_ce_{\ep}$ is a finite module over $e_{\ep}H_ce_{\ep}$, $G$ preserves finite dimensionality. Suppose $M$ is a finite dimensional simple $e_{\ep}H_ce_{\ep}$--module. Consider a composition series for $G(M)$ $$G(M) = X_0 \supset X_1 \supset \cdots \supset X_n \supset 0.$$ By \cite[Th\'{e}or\`{e}me 4.1]{dez} each finite dimensional simple $H_c$--module belongs to the category $\mathcal{O}_c$. By Lemma \ref{maps} if $L_c(\tau)$ is finite dimensional then $\tau \cong \h_i$ for some $i$. Thus Lemma \ref{faithful} shows that $F$ is exact and faithful on the finite dimensional simple modules of $\mathcal{O}_c$. We deduce a composition series $$F(G(M)) = F(X_0) \supset F(X_1)\supset \cdots \supset F(X_n) \supset 0.$$ But $F(G(M)) = e_{\ep}(H_ce_{\ep} \otimes_{e_{\ep}H_ce_{\ep}} M) = M$. Hence $F(X_1) = 0$, implying that $X_1 = 0$. Hence $G$ preserves the simplicity of finite dimensional representations. 

Let $L = G(\mathbb{C})$. By construction $L$ contains a copy of $\ep$ in degree 0. By Lemma \ref{ppsign} the unique standard module $M_c(\h_i)$ with $\ep$ appearing in degree $0$ is $M_c(\h_0)$. Thus $L$ must be a factor of $M_c(\h_0)$. In other words $L \cong L_c(\h_0)$. By Lemma \ref{altsum} we deduce that 
$$ L_c(\h_0) = \sum_{i=0}^n (-1)^i M_c(\h_i).$$
Thus we find that \begin{eqnarray*} p(L_c(\h_0),t) &=& \sum_{i=0}^n (-1)^i p(M_c(\h_i),t) \\ &=& \sum_{i=0}^n (-1)^i \frac{\binom{n}{i}t^{-N + (h+1)i}}{(1-t)^n} \\ &=& t^{-N}\frac{(1-t^{h+1})^n}{(1-t)^n}   \\ &=& t^{-N}(1+t+t^2 + \cdots + t^h)^n.\end{eqnarray*}
\end{proof}
The result of this theorem is consistent with the conjecture in \cite[Section 5]{BEG} on the character formula for representations of $H_c$ in type $A$. This conjecture has now been confirmed in \cite{BEG2}.
\section{Diagonal harmonics}
The following theorem was conjectured by Haiman for all Weyl groups, \cite[Section 7]{haiman}. For type $A$ a stronger version of the theorem was obtained by Haiman in \cite{haiman2}. Results for dihedral groups were obtained by Alfano and Reiner. 

Recall $n$ is the rank of $W$, $h$ its Coxeter number of $W$, $Q$ the root lattice associated with $W$, and $\ep$ the sign representation of $W$.
\begin{thm}
Let $W$ be a Weyl group. Let the quotient ring by diagonal invariants be $$D_W = \frac{\mathbb{C}[\h \oplus \h^*]}{\langle \mathbb{C}[\h \oplus \h^*]_+^W\rangle }.$$ Then there exists a $W$--stable quotient ring $R_W$ of $D_W$ satisfying the following properties:
\begin{enumerate} 
\item $\dim R_W = (h+1)^n$;
\item $R_W$ is $\mathbb{Z}$--graded with Hilbert series $t^{-N}(1+t+\cdots +t^h)^n$;
\item The image of $\mathbb{C}[\h]$ in $R_W$ is the classical coinvariant algebra, $\mathbb{C}[\h]^{co W}$;
\item As a $W$--module $R_W\otimes \ep$ is isomorphic to the permutation representation of $W$ on $Q/(h+1)Q$.
\end{enumerate}
\end{thm} 
\begin{proof}
Let $L= H_ce_{\ep}\otimes_{e_{\ep}H_ce_{\ep}} \mathbb{C}$ for $c=(1+h)/h$. Consider the associated graded module $gr L$. Since $gr (H_ce_{\ep}) = \mathbb{C}[\h\oplus \h^*]\ast W e_{\ep}$ and $gr (e_{\ep}H_ce_{\ep}) = e_{\ep}\mathbb{C}[\h\oplus \h^*]\ast We_{\ep}$ we have a natural surjection $$\mathbb{C}[\h\oplus \h^*]\ast We_{\ep}\otimes_{ e_{ep}\mathbb{C}[\h\oplus \h^*]\ast We_{\ep}}\mathbb{C} \longrightarrow gr L.$$ Left multiplication by $e_{\ep}$ provides a graded $\mathbb{C}[\h\oplus \h^*]\ast W$--isomorphism between $\mathbb{C}[\h\oplus \h^*]^We_{\ep}$ and $e_{\ep}\mathbb{C}[\h\oplus \h^*]\ast We_{\ep}$. Right multiplication by $e_{\ep}$ provides a graded $\mathbb{C}[\h\oplus \h^*]\ast W$--isomorphism between $\mathbb{C}[\h\oplus \h^*]\otimes \ep$ and $\mathbb{C}[\h\oplus \h^*]\ast We_{\ep}$. We deduce a graded $\mathbb{C}[\h\oplus \h^*]\ast W$--surjection $$\mathbb{C}[\h\oplus \h^*]\otimes_{\mathbb{C}[\h\oplus \h^*]^W}\mathbb{C} \longrightarrow gr L\otimes \ep.$$ We set $R_W = gr L\otimes \ep$.

By the second paragraph of the proof of Theorem \ref{hs} the element $e_{\ep}\otimes 1$ has degree $0$ in $L$. Since the grading induced by $\mathbf{h}$ gives $x\in \h^*$ degree $1$ and $y\in \h$ degree $-1$, we see that $R_W$ has the same Hilbert series as $L$, so (2) follows from Theorem \ref{hs}. Part (1) is obtained by setting $t=1$.

The image of $\mathbb{C}[\h]$ in $R_W$ corresponds to the subspace  $\mathbb{C}[\h]e_{\ep}\otimes 1$ of $L$. If $p\in \mathbb{C}[\h]^W_+e_{\ep}$ then $$p\otimes 1 = e_{\ep}pe_{\ep}\otimes 1 = e_{\ep}\otimes p.1 = 0.$$ Thus the ideal generated by $\mathbb{C}[\h]_+^W$ annihilates $e_{\ep}\otimes 1$. On the other hand, the quotient $\mathbb{C}[\h]^{co W}$ contains a unique (up to scalar) element of maximal degree $N$, say $q$, \cite[Chapter 3]{hum}. The space $\mathbb{C}q$ is the socle of $\mathbb{C}[\h]^{co W}$ since the ring of coinvariants is Frobenius (using Poincar\'{e} duality, for instance). We claim $qe_{\ep}\otimes 1\neq 0$. By (\ref{PBW}) any element of $H_c$ can be written as a sum of terms of the form $p_-wp_+$ where $p_-\in\mathbb{C}[\h^*]$, $p_+\in\mathbb{C}[\h]$ and $w\in W$. Since $p_-$ and $w$ do not increase degree, it would follow if $qe_{\ep}\otimes 1$ were zero, then $L$ could have no subspace in degree $N$. But the Hilbert series of $L$ has highest order term $t^{-N+hn} = t^N$. Thus $qe_{\ep}\otimes 1$ is non--zero and $\mathbb{C}[\h]e_{\ep}\otimes 1$ is isomorphic to $\mathbb{C}[\h]^{co W}\otimes e_{\ep}$. This proves (3).

It remains to check (4). Notice it is enough to calculate the $W$--decomposition of $L$, since passing to associated graded module is $W$--equivariant. Recall from the proof of Theorem \ref{hs} we can calculate the $W$--decomposition of $L$ from the formula $$L = \sum_{i=0}^n (-1)^iM_c(\h_i).$$ Therefore, with the obvious notation, the graded character of $w$ on $L$ is given by \begin{equation} \label{eq1} \text{ch}_L(w,t) = \sum_{i=0}^n (-1)^i \text{ch}_{M_c(\h_i)}(w,t).\end{equation}

Recall, as a graded $W$--module we have $M_c(\h_i) = \mathbb{C}[\h]\otimes \h_i$. It is well known that the graded trace of $w$ on $\mathbb{C}[\h]$ is given by $1/\det(1-tw)$. Thus \begin{equation} \label{eq2} \text{ch}_{M_c(\h_i)}(w,t) = \frac{t^{\kappa_c(\h_i)}\text{ch}_{\h_i}(w)}{\det(1-tw)}.\end{equation} Combining (\ref{eq1}), (\ref{eq2}) and the equality $\kappa_c(\h_i) = i(1+h)$ from Lemma \ref{centchar} yields $$\text{ch}_L(w,t) = t^{-N} \frac{\det(1-t^{h+1}w)}{\det (1-tw)}.$$
Evaluating this expression at $t=1$ gives the character of $L$ \begin{equation} \label{eq3} \text{ch}_L(w) = h^{\dim \ker(1-w)}.\end{equation}

On the other hand, the character of $w$ on the permutation representation $Q/(h+1)Q$ equals the number of fixed points of $w$ on $Q/(h+1)Q$. For the exceptional groups, the order of $w$ and of $h+1$ are coprime. It follows that $(Q/(h+1)Q)^w = Q^w/(h+1)Q^w$. Thus the character agrees of the permutation representation agrees with (\ref{eq3}). This proves (4) for exceptional Weyl groups.

We turn to the classical groups. The case with $W=W(A_n)$ is given in \cite[Proposition 2.5.3]{haiman}. We will prove the case $W=W(B_n)$, the case for $D_n$ being entirely similar. Recall in this case that $h=2n$, and $\h = \mathbb{C}x_1+\cdots \mathbb{C}x_n$, with $W$ acting by pemutations of the basis vectors, and multiplications by plus or minus one, \cite[Chapter 2]{hum}. Consider the polynomial ring $\mathbb{C}[\h] = \mathbb{C}[x_1,\ldots ,x_n]$. Let $I$ be the ideal generated by the vector space $V=\mathbb{C}x_1^{2n+1} + \cdots + \mathbb{C}x_n^{2n+1}$.  Observe that $V$ is isomorphic to the reflection representation  of $W(B_n)$. Since $\mathbb{C}[x_1,\ldots ,x_n]/I$ has dimension $(2n+1)^n$, it follows that $(x_1^{2n+1},\ldots ,x_n^{2n+1})$ is a homogeneous system of parameters. Therefore there is a Koszul resolution of $\mathbb{C}[\h]/I$ as a $\mathbb{C}[\h]\ast W$--module $$0\longrightarrow F_n\longrightarrow \cdots \longrightarrow F_1\longrightarrow F_0\longrightarrow \mathbb{C}[\h]/I \longrightarrow 0.$$ In this resolution $F_k$ is $\mathbb{C}[\h]\otimes \wedge^kV$. Each $F_k$ is graded by assigning degree $2n+1$ to each $x_i^{2n+1}$ and extending this to the exterior products. The maps in the resolution are homogeneous of degree $0$. We deduce that $$\mathbb{C}[\h]/I = \sum_{k=0}^n (-1)^k F_k = \sum_{i=0}^n (-1)^i M_c(\h_i).$$ Therefore it is enough to show that $\mathbb{C}[\h]/I\otimes \ep$ is isomorphic to the permutation representation on $Q/(2n+1)Q$.

The positive root vectors for $W$ are the elements $ x_i\pm x_j$ for $1\leq i<j\leq n$ and $2x_i$ for $1\leq i\leq n$. Since $2n+1$ is odd it follows that $Q/(2n+1)Q$ is isomorphic as a $W$--set to the set $S=\mathbb{Z}_{2n+1}x_1 +\cdots +\mathbb{Z}_{2n+1}x_n$. To finish the proof we will set up a bijection, $\theta$, of $W$--sets between the monomials in $\mathbb{C}[\h]/I$ and a certain $W$--stable basis of $\mathbb{C}S$. Given $s\in S$, we write $[s]$ for the corresponding element in $\mathbb{C}S$. For $1\leq i,m\leq n$ we set $$ \ep_{i,2m}= [mx_i] + [(2n+1 -m)x_i] \quad \text{ and}\quad  \ep_{i,2m+1} = [mx_i] - [(2n+1 - m)x_i].$$ Set $\ep_{i,0} = [0x_i]$. Then we define $$\theta(x_1^{m_1}\ldots x_n^{m_n}) = \sum_{i=1}^n \ep_{i,m_i}.$$ It is straightforward to check that $\theta$ has the desired properties.
\end{proof}



\begin{thebibliography}{999}
\bibitem{ariki} S.Ariki, On the decomposition numbers of the Hecke algebra of $G(m,1,n)$. J. Math. Kyoto Univ. 36 (1996), no. 4, 789--808.
\bibitem{BEG} Y.Berest, P.Etingof and V.Ginzburg, Cherednik algebras and differential operators on quasi--invariants, to appear Duke Math.J.
\bibitem{BEG2} Y.Berest, P.Etingof and V.Ginzburg, Finite dimensional representations of rational Cherednik algebras, in preparation.
\bibitem{gecketal}F.M.Bleher, M Geck and W.Kimmerle, Automorphisms of generic Iwahori--Hecke algebras and integral group rings of finite
Coxeter groups. J. Algebra 197 (1997), no. 2, 615--655.
\bibitem{bmr} M.Brou\'{e}, G.Malle and R.Rouquier, Complex reflection groups, braid groups, Hecke algebras. J. Reine Angew. Math. 500 (1998),
127--190.
\bibitem{curkil} C.W.Curtis, N.Iwahori and R.Kilmoyer, Hecke algebras and characters of parabolic type of finite groups with $(B,$ $N)$-pairs. Inst. Hautes Etudes
Sci. Publ. Math. No. 40 (1971), 81--116. 
\bibitem{dez} C.Dezelee, Representations de dimension finie de l'algebre de Cherednik rationelle, arXiv:math. RT/0111210.
\bibitem{dunkl} C.Dunkl, Differential--difference operators and monodromy representations of Hecke algebras, Pacific.J.Math. 159(2) (1993), 271--298.
\bibitem{DJO} C.Dunkl, M.F.E. de Jeu and E.Opdam, Singular polynomials for finite reflection groups, Trans. Amer.Math.Soc. 346 (1994), 237--256.
\bibitem{do} C.Dunkl and E.Opdam, Dunkl operators for complex reflection groups, preprint, arXiv:math. RT/0108185.
\bibitem{EG} P.Etingof and V.Ginzburg, Symplectic reflection algebras, Calogero--Moser space, and deformed Harish--Chandra homomorphism, Invent. Math. 147 (2002), 243--348.
\bibitem{ginz} V.Ginzburg, Principal nilpotent pairs in a semisimple Lie algebra. I, Invent.Math. 140 (2000), 511--561.
\bibitem{guay} N.Guay, Projective modules in the category $\mathcal{O}$ for the Cherednik algebra, preprint.
\bibitem{haiman} M.Haiman, Conjectures on the quotient ring by diagonal invariants, J.Algebraic Combin. 3 (1994), 17--76.
\bibitem{haiman2} M.Haiman, Vanishing theorems and character formulas for the Hilbert schemes of points in the plane, Invent. Math.149 (2002), 371--407.
\bibitem{hum} J.Humphreys, Reflection groups and Coxeter groups. 
Cambridge Studies in Advanced Mathematics, 29. Cambridge University Press, Cambridge, 1990. xii+204 pp.
\bibitem{mathas} A.Mathas, Iwahori-Hecke algebras and Schur algebras of the symmetric group. University Lecture Series, 15. American Mathematical
Society, Providence, RI, 1999. xiv+188 pp.
\bibitem{nak} H.Nakajima, Lectures on Hilbert schemes of points on surfaces. University Lecture Series, 18. American Mathematical Society, Providence,
RI, 1999. xii+132 pp.
\bibitem{sol} L. Solomon, Invariants of finite reflection groups. Nagoya Math. J. 22 1963 57--64.
\bibitem{verb} M.Verbitsky, Holomorphic symplectic geometry and orbifold singularities. Asian J. Math. 4 (2000), no. 3, 553--563. 
\bibitem{wallach} N.Wallach, Invariant differential operators on a reductive Lie algebra and Weyl group representations, J.Amer.Math.Soc. 6 (1993), 779--816.

\end{thebibliography}
\end{document}